\documentclass{amsart}

\usepackage{amssymb}
\usepackage{diagrams}
\usepackage{a4}
\usepackage{hyperref}

\let\epsilon\varepsilon

\def\MFie{i.\,e.}
\def\MFcf{cf.}
\def\MFeg{e.\,g.}
\def\ZZZ{\mathbb Z}
\def\QQQ{\mathbb Q}
\def\RRR{\mathbb R}
\def\FFF{\mathbb F}
\DeclareMathOperator{\Tor}{Tor}
\def\pair#1{{\langle #1\rangle}}
\def\MFpt{\hbox{pt}}
\def\MFRP{\mathbb{RP}}
\def\HBT{{H^*(BT)}}

\def\ttt{\mathfrak t}

\theoremstyle{plain}
\newtheorem{theorem}{Theorem}[section]
\newtheorem{proposition}[theorem]{Proposition}
\newtheorem{lemma}[theorem]{Lemma}

\theoremstyle{definition}
\newtheorem{remark}[theorem]{Remark}
\newtheorem{example}[theorem]{Example}

\numberwithin{equation}{section}

\newenvironment{acknowledgements}{\begin{trivlist}%
  \item[]{\em Acknowledgements.}\hskip0.5em}{\end{trivlist}}


\begin{document}

\title[Exact cohomology sequences]{Exact cohomology sequences\\
  with integral coefficients for torus actions}
\author{Matthias Franz and Volker Puppe}

\maketitle

\begin{abstract}
  Using methods applied by Atiyah in equivariant $K$-theory,
  Bredon obtained exact sequences for the relative cohomologies
  (with rational coefficients) of the equivariant skeletons
  of (sufficiently nice) $T$-spaces~$X$, $T=(S^1)^n$, with
  free equivariant cohomology~$H_T^*(X;\QQQ)$ over~$H_T^*(\MFpt;\QQQ)=H^*(BT;\QQQ)$.
  Here we characterise those finite $T$-CW~complexes
  with connected isotropy groups for which an analogous result
  holds with integral coefficients.
\end{abstract}

\maketitle

\section{Introduction}

\noindent
Let $T=(S^1)^n$ be an $n$-dimensional torus.
By a $T$-space we mean a finite %, connected\comment{why connected?}
$T$-CW complex, for example a smooth compact $T$-manifold.
$H^*(-)$ denotes singular cohomology with integral coefficients
unless stated otherwise.
Let~$X_i$, $-1\le i\le n$, be the equivariant $i$-skeleton of~$X$,
\MFie, the union of all orbits of dimension~$\le i$.
In particular, $X_{-1}=\emptyset$, $X_0=X^T$ and $X_n=X$.

Recall that the $T$-equivariant cohomology~$H_T^*(X)$ of~$X$ is defined
as the (singular) cohomology of the Borel construction~$X_T=(ET\times X)/T$
and that $H_T^*(X)$ is a module over~$\HBT$, where
$ET\to BT$ denotes the universal $T$-bundle.
We also consider $\ZZZ$ as an $\HBT$-module
by the augmentation~$\epsilon\colon \HBT\to\ZZZ$
that annihilates all elements of positive degrees.

The inclusion of pairs~$(X_i,X_{i-1})\hookrightarrow(X,X_{i-1})$
gives rise to a long exact sequence in equivariant cohomology
\begin{equation}\label{exact-triple}
  \cdots\to H^*_T(X, X_i)
  \to H^*_T(X, X_{i-1})
  \to H^*_T(X_i, X_{i-1})
  % \stackrel\delta
  \to H^{*+1}_T(X, X_i)
  \to\cdots,
\end{equation}
and likewise $(X_{i+1},X_{i-1})\hookrightarrow(X_{i+1},X_i)$ gives
boundary maps
\begin{equation}
  H^*_T(X_i, X_{i-1})\to H^{*+1}_T(X_{i+1}, X_i).
\end{equation}

\begin{theorem}% [Atiyah--Bredon]
  \label{main-result}
  Let $X$ be a $T$-space. %as above\comment{why `as above'}
  % with connected isotropy groups.
  Then the following four conditions are equivalent:
  \def\theenumi{\roman{enumi}}
  \begin{enumerate}
  \item \label{inclusion-fibre}% \label{first-condition}
    The inclusion of the fibre~$\iota\colon X\hookrightarrow X_T$
    induces a surjection~$\iota^*\colon H_T^*(X)\to H^*(X)$.
    Equivalently,
    the map~$H_T^*(X)\otimes_\HBT\ZZZ\to H^*(X)$ induced by~$\iota^*$
    is an isomorphism.
  \item \label{Serre-degenerates}
    The Serre spectral sequence
    of the fibration~$X\hookrightarrow X_T\to BT$
    degenerates at the $E_2$~level, \MFie, $E_2=E_\infty$.
  \item \label{Tor-j}
    $\Tor^{\HBT}_j(H_T^*(X),\ZZZ)=0$ for all~$j>0$.
  \item \label{Tor-1}
     $\Tor^{\HBT}_1(H_T^*(X),\ZZZ)=0$.
  \end{enumerate}
  The above conditions are implied by the next one.
  \begin{enumerate}
  \setcounter{enumi}4
  \item \label{Atiyah-exact}% \label{last-condition}
    The following sequence is exact:
    $$ % \label{sequence-atiyah}
      0
      \to H^*_T(X)
      \to H^*_T(X_0)
      \to H^{*+1}_T(X_1, X_0)
      \to \cdots
      \to H^{*+n}_T(X_n, X_{n-1})
      \to 0.
    $$
  \end{enumerate}
  Moreover,
  if all isotropy groups of~$X$ are connected, then this last condition
  is equivalent to the others.
\end{theorem}

The main contribution of this paper is to relate
the criteria \eqref{inclusion-fibre}--\eqref{Tor-1}
to the exact sequence given in~\eqref{Atiyah-exact}.
For field coefficients, the equivalence of \eqref{inclusion-fibre}
to~\eqref{Tor-1} is standard.

% Note that condition~\eqref{inclusion-fibre} could also be formulated
% in the following way:
% The map~$H_T^*(X)\otimes_\HBT\ZZZ\to H^*(X)$ induced by~$\iota^*$
% is an isomorphism, where the $\HBT$-module structure on~$\ZZZ$ is given
% by the augmentation~$\epsilon\colon \HBT\to\ZZZ$
% that annihilates all elements of positive degrees.

It is easy to see that for coefficients in a field~$k$ each of the
conditions \eqref{inclusion-fibre}~to~\eqref{Tor-1} above is equivalent
to $H_T^*(X;k)$ being a free $H^*(BT;k)$-module. This does not hold
for integer coefficients, as Example~\ref{example-2} below shows.

The connectedness assumption on the isotropy groups
will be used in Proposition~\ref{localization-iso}
and Lemma~\ref{finitely-generated}
and cannot be omitted in general.
But, as shown in~\cite{FranzPuppe:??},
one can allow one (finite) cyclic factor in each isotropy group
if $H_T^*(X)$ is \emph{free} over~$\HBT$.
In Section~\ref{examples} we present some examples of what can go wrong
if these conditions are not satisfied.
Moreover, one can obtain a partial result if the isotropy groups
are connected for all~$x$ in some~$X_{i}$
and not ``too much'' disconnected for~$x\not\in X_{i}$,
see Remark~\ref{weaker-result}.

\medskip

  That the freeness of~$H^*(X;\QQQ)$ implies the
  exact sequence~\eqref{Atiyah-exact} (with rational coefficients)
  was pointed out by Bredon~\cite{Bredon:74}, based on an analogous result
  in equivariant $K$-theory due to Atiyah \cite[Lecture~7]{Atiyah:74}.
  We therefore refer to sequence~\eqref{Atiyah-exact} as
  the ``Atiyah--Bredon sequence''.

  Atiyah's proof relies on the theory of Cohen--Macaulay modules.
  In~\cite{FranzPuppe:??} the same method is used
  to give a version for other coefficients, which also allows
  for non-connected isotropy groups under appropriate assumptions
  (see the discussion preceding Example~\ref{example-4} below).
  In the present paper, we instead apply fairly standard homological algebra
  which makes essential use of the grading and gives a more complete
  picture in the case of connected isotropy groups and integer coefficients.
  This approach was inspired by
  Barthel--Brasselet--Fieseler--Kaup, who prove a result analogous to the
  rational version of Theorem~\ref{main-result}
  in the context of (intersection) cohomology for fans
  \cite[Theorem~4.3]{BarthelBrasseletFieselerKaup:02}.

\medskip
% \begin{remark}\rm
  For $H_T^*(X;\QQQ)$ a free $H^*(BT;\QQQ)$-module, the exactness of the
  Atiyah--Bredon sequence at~$H_T^*(X;\QQQ)$
  is an immediate consequence of the Localisation Theorem, and the
  exactness at~$H_T^*(X_0;\QQQ)$ follows directly from the so-called
  Chang--Skjelbred lemma~%
  \cite[Lemma~2.3 \& Proposition~2.4]{ChangSkjelbred:74}.
  Goresky--Kottwitz--MacPherson~\cite{GoreskyKottwitzMacPherson:98}
  and others made very effective use of the Chang--Skjelbred result
  in calculating the (equivariant) cohomology of certain $T$-spaces,
  sometimes even with coefficients in~$\ZZZ$.
% \end{remark}

\begin{remark}\label{comment-real}\rm
  The proof for Theorem~\ref{main-result} given below simplifies
  if one takes rational (or real) coefficients instead of~$\ZZZ$.
  In this case the assumption about the connectedness of the isotropy groups
  can be dropped completely. The reason is that for any isotropy group~$T_x$
  one has that $H^*(BT_x;\QQQ)$ is naturally isomorphic to~$H^*(BT_x^0;\QQQ)$
  where $T_x^0\subset T_x$ denotes the identity component.
\end{remark}

\begin{remark}\rm
  As with the Atiyah--Bredon method (\MFcf~\cite{Bredon:74}),
  our arguments can also be used to prove
  similar results for $p$-tori and coefficients in~$\FFF_p$.
  Again, the proof % is analogous to the one given below, but
  simplifies since we can take as extension of~$\FFF_p$ any infinite field
  of characteristic~$p$ (\MFeg, the algebraic closure~$\bar\FFF_p$)
  instead of the less familiar extension of~$\ZZZ$ introduced
  in Section~\ref{new-ground-ring}.
  We nevertheless concentrate on the (compact) torus case
  and integer coefficients since we get new results in this context.
  Variants for compact tori and coefficients in~$\FFF_p$ or
  in a subring of~$\QQQ$ can be found in~\cite{FranzPuppe:??}.
\end{remark}

\begin{remark}\rm
  The assumption that $X$ is a finite $T$-CW~complex could be weakened
  if one is willing to apply more technical machinery along the lines
  of~\cite[Chapter~3]{AlldayPuppe:93}.
\end{remark}

\medbreak

\begin{acknowledgements}
  The authors thank Walter Baur, Karl-Heinz Fieseler
  % and Jean-Claude Hausmann
  and Sue Tolman
  for helpful discussions, and the referees for their comments
  and suggestions.
\end{acknowledgements}

\section{Some homological algebra}

\noindent
In this section we collect a few more or less known results from homological
algebra for later use.

We work over an arbitrary ground ring (commutative, with unit)~$k$.
Unless stated otherwise, tensor products and direct sums are taken over~$k$.
We set $A=k[t_1,\ldots,t_r]$ where each $t_i$ has degree~$2$, and we call
elements in~$A^2$ (= the $k$-module generated by~$t_1$,~\ldots,~$t_r$) linear.
All complexes and (graded) modules are assumed to be bounded from below.

The following standard fact is immediate from the long exact $\Tor$~sequence.

\begin{lemma}\label{Tor-triple}
  Let $0\to L\to M\to N\to0$ be a short exact sequence of~$A$-modules
  such that $\Tor^A_j(M,k)=0$ for all~$j>i$. Then
  $\Tor^A_{j+1}(N,k)=\Tor^A_j(L,k)$ for all~$j>i$.
\end{lemma}

An important role in this section is played
by the short exact sequence of $A[t]$-modules
\begin{equation}
  \label{kt-kt-k}
  0\to k[t]\stackrel{\cdot t}\longrightarrow k[t]\to k\to 0.
\end{equation}

\begin{lemma}\label{Tor-vanishing}
  Let $M$ be an $A[t]$-module (where $t$ is an indeterminate of degree~$2$),
  considered as an $A$-module by restriction, $A\hookrightarrow A[t]$.
  Then $\Tor^A_1(M,k)$ vanishes if one of the following conditions holds:
  \def\theenumi{\roman{enumi}}
  \begin{enumerate}
  \item \label{Tor-vanishing-1}
    $\Tor^{A[t]}_1(M,k)=0$,
  \item \label{Tor-vanishing-2}
    $\Tor^{A[t]}_2(M,k)=0$ and $M$ is finitely generated over~$A$.
  \end{enumerate}
\end{lemma}

\begin{proof}
  We first remark that $\Tor^A_j(M,k)=\Tor^{A[t]}_j(M,k[t])$.
  This follows from the fact that any free resolution of~$k$ over~$A$,
  tensored over~$k$ with~$k[t]$, gives a free resolution of~$k[t]$
  over~$A[t]$.

  The sequence~\eqref{kt-kt-k} induces the exact sequence 
  $$
    \cdots\to\Tor^{A[t]}_2(M,k)
      \to\Tor^{A[t]}_1(M,k[t])
      \stackrel{\theta}\longrightarrow\Tor^{A[t]}_1(M,k[t])
      \to\Tor^{A[t]}_1(M,k)\to\cdots.
  $$

  If $\Tor^{A[t]}_1(M,k)=0$, then the map~$\theta$ is surjective.
  But it is also of positive degree, hence zero
  since all modules are bounded from below.
  Therefore, $\Tor^A_1(M,k)=\Tor^{A[t]}_1(M,k[t])=0$.

  If $\Tor^{A[t]}_2(M,k)=0$, then $\theta$ is injective.
  But $\Tor^{A[t]}_1(M,k[t])=\Tor^{A}_1(M,k)$
  is a finitely generated $A$-module such that
  the $A$-action factors through the augmentation~$A\to k$.
  Hence $\Tor^A_1(M,k)$ is a finitely generated $k$-module,
  which is necessarily bounded from above. Therefore, the map~$\theta$,
  which raises degrees, cannot be injective unless $\Tor^A_1(M,k)$
  vanishes.
\end{proof}

The equivalence of conditions \eqref{Tor-j}~and~\eqref{Tor-1}
of Theorem~\ref{main-result}
% will somehow drop out of our proof
% in Section~\ref{proof-main-result}. But since it
is purely algebraic in nature
% we give a short direct proof of it here.
and moreover well-known if $k$ is a field.
% If $k$ is a field, this result is well-known.
The conditions are then equivalent to $M=H^*_T(X)$~being free
over~$A=\HBT$. This does not hold in general,
\MFcf~Remark~\ref{extension-free} and Example~\ref{example-2}.

\begin{proposition}\label{higher-Tor-vanish}
  Let $M$ be an $A$-module.
  %MF060116
  % Then $\Tor^A_1(M,k)=0$ if and only if $\Tor^A_j(M,k)=0$ for all~$j>0$.
  Then $\Tor^A_1(M,k)$ vanishes if and only if $\Tor^A_j(M,k)$ vanishes for all~$j>0$.
\end{proposition}

\begin{proof}
  We proceed by induction on the number~$r$ of indeterminates
  of~$A=k[t_1,\ldots,t_r]$. For~$r\le1$ there is nothing to prove.

  So assume that we have shown the equivalence for all $A$-modules,
  and let $M$ be an $A[t]$-module such that $\Tor^{A[t]}_1(M,k)=0$.
  Again, we take a portion of the long exact $\Tor$~sequence
  induced by the sequence~\eqref{kt-kt-k},
  $$
    \cdots\to\Tor^A_2(M,k)
      \to\Tor^{A[t]}_2(M,k)
      \to\Tor^A_1(M,k)
      \to\Tor^A_1(M,k)
      \to\Tor^{A[t]}_1(M,k)=0.
  $$
  (Recall that $\Tor^{A[t]}_1(M,k[t])=\Tor^A_1(M,k)$.)
  As in the proof of Lemma~\ref{Tor-vanishing}\,(\ref{Tor-vanishing-1}), we see that
  $\Tor^A_1(M,k)$ vanishes, hence by assumption
  $\Tor^A_j(M,k)$ for~$j>0$. This implies that $\Tor^{A[t]}_j(M,k)$
  vanishes as well for~$j>0$.
\end{proof}

\begin{remark}\label{extension-free}\rm
  Since a free resolution of a $k$-module~$N$ extends to a free
  resolution of the extended module~$N\otimes A$, one gets
  $\Tor^A_j(N\otimes A,k)=0$ for~$j>0$. (The converse does not hold as
  can be seen from Example~\ref{example-2}.)

  More generally, consider the subalgebras~$A'=k[t_1,\ldots,t_i]$
  and $A''=k[t_{i+1},\ldots,t_r]$ of~$A$. % Then $A=A'\otimes A''$.
  Let $M'$ be an $A'$-module. Since the Koszul resolution of~$k$
  over~$A'$ has length~$i$, $\Tor^{A'}_j(M',k)$
  vanishes for~$j>i$.
  Tensoring a free resolution of the $A'$-module~$M'$
  with~$A''$ gives a free resolution of~$M'\otimes A''$
  over~$A=A'\otimes A''$.
  By the K\"unneth theorem, this implies
  $\Tor^A_j(M'\otimes A'',k)=\Tor^{A'}_j(M',k)\otimes A''=0$
  for~$j>i$.
\end{remark}

\begin{lemma}\label{localisation-injective}
  Let $M$ be an $A$-module.
  % If $\Tor^A_1(M,k)=0$,
  If $\Tor^A_1(M,k)$ vanishes,
  then the localisation map~$M\to S^{-1}M$
  is injective for the multiplicative subset~$S\subset A$
  generated by the linear elements which can be extended to a basis
  of~$A^2$.
\end{lemma}

\begin{proof}
  Assume that the localisation map is not injective. Then there exist
  a non-zero~$m\in M$ and a linear element~$t\in S$ with~$t m=0$.
  % We choose such a pair with $m$ of minimal degree.
  Consider again the exact sequence~\eqref{kt-kt-k}.
  Tensoring with~$M$ over~$k[t]$ gives a map~$M\to M$ which sends $m'$
  to~$t m'$. This map is not injective because $t m=0$.
  Hence $\Tor^{k[t]}_1(M,k)\ne0$. Since we may extend $t=t_1$
  to a basis~$(t_1,\ldots,t_n)$ of~$A^2$, we get
  by Lemma~\ref{Tor-vanishing}\,(\ref{Tor-vanishing-1})
  $$
    \Tor^{k[t_1]}_1(M,k)\ne0
      \;\Rightarrow\;\Tor^{k[t_1,t_2]}_1(M,k)\ne0
      \;\Rightarrow\;\cdots
      \;\Rightarrow\;\Tor^{k[t_1,\ldots,t_r]}_1(M,k)\ne0,
  $$
  which contradicts the assumption~$\Tor^A_1(M,k)=0$.
\end{proof}

\section{A new ground ring}
\label{new-ground-ring}

\noindent
In the course of our proof we want to choose a common complement
for a finite set of direct summands of the same rank in~$k^n$.
This is no problem if $k$ is $\QQQ$ or~$\RRR$ (or any other infinite field),
but it is not possible in general for~$k=\ZZZ$. Therefore, we extend
the coefficients, a step which would be unnecessary in the case
of an infinite field.

Note that in Theorem~\ref{main-result}
we may replace $\ZZZ$ by any $\ZZZ$-algebra~$k$
which is a flat and faithful $\ZZZ$-module, \MFie,
such that the functor~$-\otimes_\ZZZ k$ is exact and faithful.
In other words, $k$ is $\ZZZ$-torsion-free,
and for all
% (finitely generated)
$\ZZZ$-modules~$N$
one has $N\otimes_\ZZZ k=0$ if and only if $N=0$.
Indeed, for such a~$k$ we have by the Universal Coefficient Theorem:
\begin{enumerate}
\item $H^*(C\otimes_\ZZZ k)=H^*(C)\otimes_\ZZZ k$
  for $C$, a complex,
% \item $M\otimes_\ZZZ k=0$ if and only if $M=0$,
\item
  $\Tor^{H^*(BT;k)}(M\otimes_\ZZZ k,k)=\Tor^{H^*(BT;\ZZZ)}(M,\ZZZ)\otimes_\ZZZ k$
  for $M$, an $H^*(BT;\ZZZ)$-module.
\end{enumerate}
Hence all conditions of Theorem~\ref{main-result} hold
over~$\ZZZ$ if and only if they hold over~$k$.

We will use this
to pass from $\ZZZ$ to the localisation of% the ungraded algebra
~$\ZZZ[s_1,\ldots,s_n]$
with respect to all homogeneous elements % in~$\ZZZ[s_1,\ldots,s_n]$
which are \emph{$\ZZZ$-reduced}, \MFie, which are
% Z-reduced earlier!!!
not divisible by any integer~$>1$.
We will call this new (ungraded) algebra~$k$.
Since $\ZZZ[s_1,\ldots,s_n]$ is factorial,
no prime integer in~$k$ is invertible. 
This implies in particular that
$N\otimes_\ZZZ k\neq0$
% the canonical map~$N\to N\otimes_\ZZZ k$ is injective
for~$N=\ZZZ/m\ZZZ$, hence for all% (finitely generated)
~$N\neq0$.
Because $k$ is $\ZZZ$-torsion-free as well, it is flat and faithful over~$\ZZZ$.

% This new algebra (without grading), which we will call $k$,
% has no $\ZZZ$-torsion and no prime integer is invertible.
% (Remember that $\HBT\cong\ZZZ[s_1,\ldots,s_n]$ is factorial.)
% This implies in particular that the canonical map~$M\to M\otimes_\ZZZ k$
% is injective for~$M=\ZZZ/m\ZZZ$, hence for all (finitely generated)~$M$.

\medskip

Let $W_{n-i}\subset k^n$ be the submodule generated by
$$
  w_1=(s_1,\ldots,s_n),\: w_2=(s_1^2,\ldots,s_n^2),\: \ldots,\:
  w_i=(s_1^i,\ldots,s_n^i).
$$
It has the following nice property:

\begin{lemma}\label{extension}
  Let $V\subset\ZZZ^n$ a direct summand
  and denote the induced direct summand of~$k^n$ by~$V'$.
  \begin{enumerate}
  \item \label{extension-1}
    If $\dim V= i$, then $V'\oplus W_i=k^n$.
  \item \label{extension-2}
    If $\dim V>i$,
    then $V'\cap W_i$ contains an element which can be extended
    to a basis of~$W_i$.
  \end{enumerate}
\end{lemma}

\begin{proof}$ $
% \begin{enumerate}
% \item 
\eqref{extension-1}
  Choose a basis of~$V$ and
  consider the % $(n-i)\times n$-
  matrix~$U$ whose rows are these basis vectors.
  Then the greatest common divisor of all $i\times i$-minors of~$U$
  is~$1$. %\comment{REF?}
  Let $U'$ be the $n\times n$-matrix obtained by appending
  the row vectors~$w_1$,~\ldots,~$w_{n-i}$ to~$U$.
  The determinant of~$U'$
  % completed by the vectors~$w_1$,~\ldots,~$w_{n-i}$,
  is a $\ZZZ$-reduced element of~$\ZZZ[s_1,\ldots,s_n]$,
  hence invertible in~$k$.
  This proves the claim.

% \item 
\eqref{extension-2}
  We may assume $\dim V=i+1$. By part~(\ref{extension-1}),
  we have $V'\oplus W_{i+1}=k^n$.
  Therefore, $w_{n-i}$ can be written as the sum of some~$v'\in V'$
  and a linear combination~$w$ of~$w_1$,~\ldots,~$w_{n-i-1}$.
  Hence $v'=w_{n-i}-w\in V'\cap W_i$
  can be extended to the basis~$(w_1,\ldots,w_{n-i-1},v')$ of~$W_i$.
% \end{enumerate}
\end{proof}

From now on to the end of the proof of Theorem~\ref{main-result}
in Section~\ref{proof-main-result},
$H^*(\cdot)$ denotes cohomology with coefficients
in the new ground ring~$k$.
We set $A=\HBT=k[t_1,\ldots,t_n]$.
The $k$-module~$A^2$ is identified with~$k^n$,
using the basis~$(t_1,\ldots,t_n)$.
In particular, we consider each~$W_i$ as contained in~$A^2$.

\medskip

We will need the following refined version
of the Localisation Theorem.

\begin{proposition}\label{localization-iso}
  Let $X$ be a $T$-space as above with connected isotropy groups.
  Then the map
  $$
    H_T^*(X,X_{i-1})\to H_T^*(X_i,X_{i-1})
  $$
  becomes an isomorphism after localisation
  with respect to the multiplicative subset~$S_i$ generated
  by all linear elements in~$W_i$ which can be extended
  to a basis of~$A^2\cong k^n$.
\end{proposition}

\begin{proof}
  For~$x\in X$ let $V_x$ be the kernel of the map~$H^2(BT)\to H^2(BT_x)$,
  induced by the inclusion~$T_x\hookrightarrow T$. Since $T_x$ is assumed
  to be connected, $V_x$ is a direct summand of~$H^2(BT)$
  with~$\dim V_x+\dim T_x=n$.
  By a standard argument (induction over attaching equivariant cells
  to~$X_i$ in order to get $X$, \MFcf~\cite[III.3.8]{tomDieck:87}
  or \cite[3.1.6]{AlldayPuppe:93}),
  it suffices to show that $S_i^{-1}H_T^*(T/T_x)=0$
  for all~$x\in X\setminus X_i$.
  Since $\dim V_x>i$ for such~$x$, Lemma~\ref{extension}\,(\ref{extension-2}) gives
  us an~$s\in S_i\cap V_x$ that annihilates $H_T^*(T/T_x)=H^*(BT_x)$.
\end{proof}

\begin{remark}\label{weaker-1}\rm
  It can be seen directly from the proof of Proposition~\ref{localization-iso}
  that, for a fixed~$i=0$,~\ldots,~$n$, instead of demanding the connectedness
  of all isotropy groups it suffices that
  for~$x\in X\setminus X_i$ the kernel of~$H^2(BT)\to H^2(BT_x)$ contains
  a direct summand of dimension~$>i$. This means that
  for all~$x\in X\setminus X_i$ the isotropy group~$T_x$ is contained
  in a subtorus of dimension~$n-(i+1)$, \MFeg, for~$i=0$ all isotropy groups
  different from~$T$ should be contained in a proper subtorus.
\end{remark}

Let $A_i\subset A$ be the subalgebra generated by~$W_i\subset A^2$.
Note that $A_i$ is noetherian because $k$ is.

\begin{lemma}\label{finitely-generated}
  Suppose that all isotropy groups of~$X$ are connected.
  Then the $A$-module $H_T^*(X, X_{i-1})$, considered as $A_i$-module
  via restriction, is finitely generated.
\end{lemma}

\begin{proof}
  Since $A_i$ is noetherian, it suffices
  (using long exact cohomology sequences) to show that,
  for each~$x\in X\setminus X_{i-1}$, $H_T^*(T/T_x)=H^*(BT_x)$
  is finitely generated over~$A_i$.
  But since $T_x$ is assumed to be connected,
  $\dim V_x\ge i$ and $V_x+W_i=A^2$ for all~$x\in X\setminus X_{i-1}$.
  So $H^*(BT_x)$ is a quotient of~$A_i$
  because the map~$W_i\to H^2(BT)\to H^2(BT_x)$ is surjective.
\end{proof}

\begin{remark}\label{weaker-2}\rm
  Again, for a fixed~$i=0$,~\ldots,~$n$, one can see directly from the proof
  of Lemma~\ref{finitely-generated} that it holds under a weaker condition
  than the connectedness of the isotropy groups. One needs that $V_x+W_i=A^2$
  for all~$x\in X\setminus X_{i-1}$. This holds if $T_x$ is contained in
  a subtorus of dimension~$n-i$ for all~$x\in X\setminus X_{i-1}$.
\end{remark}

\goodbreak

\begin{lemma}\label{Tor-dim-Xi}
  The following torsion products vanish:
  \begin{enumerate}
  \item\label{Tor-dim-Xi-1}
    $\Tor^A_j(H_T^*(X_i,X_{i-1}),k)=0$ for~$j>i$.
  \item\label{Tor-dim-Xi-2}
    $\Tor^{A_i}_j(H_T^*(X_i,X_{i-1}),k)=0$ for~$j>0$.
  \end{enumerate}
\end{lemma}

\begin{proof}
  The quotient~$X_i/X_{i-1}$ is a bouquet of spaces, % with base points,
  each of them fixed by a subtorus of codimension~$i$.
  Let $Y$ % $(Y,\MFpt)$
  be such a space, % with base point,
  fixed by~$T''\subset T$, and let $T'$ be a torus complement to~$T''$ in~$T$.
  Since $\HBT=H^*(BT')\otimes H^*(BT'')$
  and $H_T^*(Y,\MFpt)=H_{T'}^*(Y,\MFpt)\otimes H^*(BT'')$,
  we are in the situation of Remark~\ref{extension-free}
  with $A'=H^*(BT')$, $A''=H^*(BT'')$,
  $M'=H_{T'}^*(Y,\MFpt)$ and $M'\otimes A''=H_T^*(Y,\MFpt)$.
  Hence, $\Tor^A_j(H_T^*(Y,\MFpt),k)=0$ for~$j>i$.
  Because $H_T^*(X_i,X_{i-1})$ is a direct sum of terms
  of the form~$H_T^*(Y,\MFpt)$, this implies (\ref{Tor-dim-Xi-1}).

%   For the second part, we use the extension of coefficients from~$\ZZZ$ to~$k$
%   in an essential way.
  In order to prove~(\ref{Tor-dim-Xi-2}), it suffices to show that
  each~$H_T^*(Y,\MFpt)=M'\otimes A''$ is an extended $A_i$-module,
  again by Remark~\ref{extension-free}.
  One has $H^2(BT)=H^2(BT')\oplus W_i$
  by Lemma~\ref{extension}\,(\ref{extension-1}), hence $A=A'\otimes A_i$.
%   Since $A=A'\otimes A''$
  Therefore, any $A'$-equivariant map from~$M'$
  to some $A$-module~$N$ uniquely extends to an $A$-equivariant
  map~$M'\otimes A_i\to N$. This applies in particular to the
  canonical inclusion~$M'\to M'\otimes A''$.
  Since we also have $A=A'\otimes A''$, we can reverse the roles
  of $A_i$~and~$A''$ to get a map in the other direction.
  By uniqueness, the compositions must be the identities
  on $M'\otimes A_i$ and $M'\otimes A''$, respectively.
  Hence, both $A$-modules are isomorphic.
%   One has $H^2(BT)=H^2(BT')\oplus W_i$
%   by Lemma~\ref{extension}\,(\ref{extension-1}). % hence $A=A'\otimes A_i$.
%   Since $W_i$ and $H^2(BT'')$ are both complements to~$H^2(BT')$,
%   the composition $W_i\to H^2(BT)\to H^2(BT'')$ is an isomorphism,
%   hence also the induced map of
%   algebras~$f\colon A_i\to A\to A''=A/\langle(A')^{>0}\rangle$.
%   Consider the map of $A$-modules
%   $$
%     g\colon M\otimes A_i\to M\otimes A'',
%     \quad
%     m\otimes a\mapsto a\cdot(m\otimes1)
%   $$
%   where the dot denotes the restriction of the $A$-action
%   on~$M\otimes A''$ to~$A_i$.
%   If we filter both sides by degree with respect to~$A_i$ and $A''$,
%   then the induced map between the associated graded modules
%   is $1_M\otimes f$. The latter map is bijective, hence so is~$g$.
\end{proof}

\section{Proof of Theorem~\protect\ref{main-result}}
\label{proof-main-result}

\def\implication#1#2#3{$\hbox{\eqref{#1}}#2\hbox{\eqref{#3}}$}
\def\proofsection#1#2#3{\noindent\implication{#1}{#2}{#3}:}

\noindent
The equivalence~$\hbox{\eqref{Tor-j}}\Leftrightarrow\hbox{\eqref{Tor-1}}$
is a special case of Proposition~\ref{higher-Tor-vanish}.
We will show the implications
$\hbox{\eqref{inclusion-fibre}}
\Rightarrow\hbox{\eqref{Serre-degenerates}}
\Rightarrow\hbox{\eqref{Tor-j}}
\Rightarrow\hbox{\eqref{inclusion-fibre}}$,
$\hbox{\eqref{Atiyah-exact}}
\Rightarrow\hbox{\eqref{Tor-j}}$
and, assuming that all isotropy groups are connected,
$\hbox{\eqref{Tor-1}}\Rightarrow\hbox{\eqref{Atiyah-exact}}$.
The steps of the first chain of implications are either standard or quite easy.
It is only in the most involved part of the proof,
$\hbox{\eqref{Tor-1}}\Rightarrow\hbox{\eqref{Atiyah-exact}}$,
% that we assume all isotropy groups to be connected and
that we need the extension of the coefficients from~$\ZZZ$ to~$k$.
Recall that $A=\HBT$ (with coefficients in~$k$).

\medskip

\proofsection{inclusion-fibre}\Rightarrow{Serre-degenerates}
% This is standard. % and could as well be shown with $\ZZZ$~coefficients.
% The restriction to the fibre~$\iota\colon H_T^*(X)\to H^*(X)$ factors
% through~$H_T^*(X)\otimes_A k$, in fact the edge homomorphism in the
% Serre spectral sequence of the Borel construction is given by
% \let\inj\to
% \let\surj\to
% $$
%   H_T^q(X)\surj H_T^q(X)\otimes_A k=E_\infty^{0,q}\inj\cdots\inj E_r^{0,q}
%   \inj\cdots\inj E_2^{0,q}=H^q(X).
% $$
% So $E_2=E_\infty$ implies that $\iota^*$ is surjective, \MFie,
% $\hbox{\eqref{Serre-degenerates}}\Rightarrow\hbox{\eqref{inclusion-fibre}}$.
The surjectivity of~$\iota^*$ implies that all boundaries
starting at~$E_r^{0,*}$ must vanish. Since $E_2^{p,q}=A^p\otimes H^q(X)$,
% the derivation property of the boundary
the $A$-linearity of the boundary maps
gives $d_r^{*,*}=0$ for all~$r\ge2$,
\MFie, $E_2=E_\infty$. %, \MFie,
% $\hbox{\eqref{inclusion-fibre}}\Rightarrow\hbox{\eqref{Serre-degenerates}}$.

\medskip

\proofsection{Serre-degenerates}\Rightarrow{Tor-j}
{\def\FFF{\mathcal F}
The abutment~$E_\infty$ of the Serre spectral sequence
can be considered as the (bi)graded $A$-module associated to
a filtration of~$H_T^*(X)$ by $A$-submodules~$\FFF_q$
such that $E_\infty^{*,q}=\FFF_q/\FFF_{q-1}$,
\MFcf~\cite[bottom of p.~554]{Quillen:71}.
Since $X$ is assumed to be a finite $T$-CW~complex,
this filtration is finite.
For~$E_\infty=E_2=A\otimes H^*(X)$ we have
$\Tor^A_j(E_\infty^{*,q},k)=0$ for~$j>0$ and any~$q$
(see Remark~\ref{extension-free}).
% Since $E_2=E_\infty$, $E_2$ is the associated graded $A$-module
% of a finite filtration on~$H_T^*(X)$, which starts
% with the extended $A$-module~$\FFF_0=A\otimes H^0(X)$.
Hence, via the long exact $\Tor$~sequences coming from
the short exact sequences $0\to\FFF_{q-1}\to\FFF_q\to E_\infty^{*,q}\to0$
one gets by induction that $\Tor^A_j(H_T^*(X),k)=0$ for~$j>0$.}

\medskip

\proofsection{Tor-j}\Rightarrow{inclusion-fibre}
This follows immediately from the
Eilenberg--Moore spectral sequence
$\Tor^A_*(H_T^*(X),k)\Rightarrow H^*(X)$.
A version of the Eilenberg--Moore theorem suitable for our purpose
% with integer coefficients
can be found in~\cite[Corollary~3.5]{GugenheimMay:74}.
% Alternatively, one can use
% the ``singular Cartan model''~$C^*(X)\tildeotimes \HBT$ for~$H_T^*(X)$,
% where ``$\tildeotimes$'' indicates a twisted differential. In this model,
% the restriction to the fibre~$H_T^*(X)\to H^*(X)$
% is induced by the canonical map~$C^*(X)\tildeotimes \HBT\to C^*(X)$.
% The latter factors as
% $$
%   C^*(X)\tildeotimes \HBT\to C^*(X)\tildeotimes K^*\to C^*(X),
% $$
% where $K^*=\HBT\tildeotimes H^*(T)\to k$ is
% the Koszul resolution of~$k$ over~$\HBT$.
% Using suitable filtrations on the complexes and condition~\eqref{Tor-j},
% one sees that in cohomology one gets the maps
% $$
%   H_T^*(X)
%   \to \Tor^{\HBT}_*(H_T^*(X),k)=H_T^*(X)\otimes_{\HBT}k
%   = H^*(X).
% $$
% \textbf{[This is less than optimal, even if references were added.]}

% \medskip

% \proofsection{Tor-j}\Rightarrow{Tor-1}
% is trivial.

\medskip

Before entering the central part of our proof, we reformulate
condition~\eqref{Atiyah-exact}.

\begin{lemma}
  The Atiyah--Bredon sequence is exact
  if and only if for all~$0\le i\le n$
  the long exact sequence~\eqref{exact-triple}
  splits into short exact sequences 
  \begin{equation}\label{short-sequence}
    0
    \longrightarrow H^*_T(X,X_{i-1})
    \stackrel{\alpha_i}\longrightarrow H^*_T(X_i,X_{i-1})
    \stackrel{\delta_i}\longrightarrow H^{*+1}_T(X,X_i)
    \longrightarrow 0.
  \end{equation}
\end{lemma}

\begin{proof}
  The idea is to combine the Atiyah--Bredon sequence
  as follows
  with the short sequences~\eqref{short-sequence}
  to a commutative diagram:
%   This is proved by chasing through a diagram which is obtained
%   by splicing the Atiyah--Bredon sequence
%   into short sequences %~\eqref{short-sequence}
%   and using the fact that the latter are part
%   of the long exact sequence~\eqref{exact-triple}.
%   The beginning of the diagram looks like this:
  {\def\zero{\mkern1mu}%
  \let\ddots=0
  \let\iddots=0
  \begin{diagram}[size=1.65em,objectstyle=\scriptstyle]
    & & & & \ddots & & & & \iddots & & & & \ddots & & & & \iddots\\
    & & & & & \rdTo & & \ruTo & & & & & & \rdTo & & \ruTo \\
    & & & & & & H_T^{*+1}(X,X_0) & & & & & & & & H_T^{*+3}(X,X_2)\\
    & & & & & \ruTo & & \rdTo & & & & & & \ruTo & & \rdTo\\
    0 & \rTo & H_T^*(X) & \rTo & H_T^*(X_0) & & \rTo & & H_T^{*+1}(X_1,X_0) & & \rTo & & H_T^{*+2}(X_2,X_1) & & \rTo & & \cdots  \\
    & & & \ruTo & & & & & & \rdTo & & \ruTo \\
    & & H_T^*(X) & & & & & & & & H_T^{*+2}(X,X_1) \\
    & \ruTo & & & & & & & & \ruTo & & \rdTo \\
    \iddots & & & & & & & & \iddots & & & & \ddots \\
  \end{diagram}}
  \noindent
  Assume that the Atiyah--Bredon sequence is exact.
  Using that the diagonal sequences come from
  the long exact sequence~\eqref{exact-triple},
  one shows inductively that they are in fact short exact.
  The reverse direction is done by diagram chase.
\end{proof}

\proofsection{Atiyah-exact}\Rightarrow{Tor-j}
By downward induction on~$i$, we prove the statement
\begin{equation}\label{Tor-i-j}
  \Tor^A_j(H_T^*(X, X_{i-1}),k)=0
  \quad
  \hbox{for all~$j>i$.}
\end{equation}
For~$i=n+1$ there is nothing to prove.
The induction step follows from  Lemma~\ref{Tor-dim-Xi}\,(\ref{Tor-dim-Xi-1})
and the long exact $\Tor$~sequence
applied to the short exact sequence~\eqref{short-sequence}.
Since $X_{-1}=\emptyset$, the case~$i=0$ is condition~\eqref{Tor-j}.
% We apply the long exact $\Tor$~sequence
% to the short exact sequence~\eqref{short-sequence} and
% use Lemma~\ref{Tor-dim-Xi}. Assuming~\eqref{Tor-i-j},
% $\Tor^A_j(H_T^*(X,X_{i-1}),k)=0$ for~$j>i$
% (which holds for~$i=n$ by Lemma~\ref{Tor-dim-Xi}),
% the long exact $\Tor$~sequence for~\eqref{short-sequence}
% gives $\Tor^A_j(H_T^*(X,X_{i-2}),k)=0$ for~$j>i-1$.
% Since $X_{-1}=\emptyset$, one gets for~$i=0$
% that $\Tor^A_j(H_T^*(X),k)=0$ for~$j>0$.

\medskip

\proofsection{Tor-1}\Rightarrow{Atiyah-exact}
Here we assume all isotropy groups to be connected.
We will show that for~$0\le i\le n$ the condition
\begin{equation}\label{Tor-i}
  \Tor^{A_i}_1(H_T^*(X,X_{i-1}),k)=0
\end{equation}
(which is satisfied for~$i=0$ by hypothesis) implies that
the sequence~\eqref{short-sequence} is exact and that
% the exactness of~\eqref{short-sequence} implies that
\eqref{Tor-i} holds for~$i+1$ instead of~$i$.
This would prove the claim by induction.

Suppose that \eqref{Tor-i} is true.
By Lemma~\ref{localisation-injective},
the localisation map~$H_T^*(X,X_{i-1})\to S_i^{-1}H_T^*(X,X_{i-1})$
is injective.
Moreover, Proposition~\ref{localization-iso} implies
that the localised map $S_i^{-1}\alpha_i$ is an isomorphism.
Hence $\alpha_i$ is injective, and \eqref{short-sequence} is exact.

Applying $\Tor^{A_i}_*(-,k)$ to~\eqref{short-sequence}
gives the exact sequence
$$
  \Tor^{A_i}_2(H_T^*(X_i,X_{i-1}),k)
  \longrightarrow\Tor^{A_i}_2(H_T^*(X,X_i),k)
  \longrightarrow\Tor^{A_i}_1(H_T^*(X,X_{i-1}),k).
$$
The first term vanishes by  Lemma~\ref{Tor-dim-Xi}\,(\ref{Tor-dim-Xi-2})
and the third one by hypothesis.
Therefore, $\Tor^{A_i}_2(H_T^*(X,X_i),k)$ vanishes as well,
as does $\Tor^{A_{i+1}}_1(H_T^*(X,X_i),k)$
by Lemmas~\ref{finitely-generated}
and~\ref{Tor-vanishing}\,(\ref{Tor-vanishing-2})
since $A_{i+1}=A_i[t]$ for some~$t\in A^2$.
% % This is the most involved part of the proof. In fact it is only here
% % that we need the extension of the coefficients from~$\ZZZ$ to~$k$.
% We will show by induction that the sequences~\eqref{short-sequence}
% are exact. By assumption $\Tor^A_1(H_T^*(X),k)=0$, and hence
% the localisation map~$H_T^*(X)\to S^{-1}H_T^*(X)$ is injective
% by Lemma~\ref{localisation-injective}.
% By the Localisation Theorem (\MFcf~Proposition~\ref{localization-iso})
% $S^{-1}H_T^*(X)=S^{-1}H_T^*(X_0)$, and therefore $\alpha_0$
% in~\eqref{short-sequence} must be injective, and \eqref{short-sequence}
% is exact for~$i=0$.
% Lemma~\ref{Tor-triple} now implies that $\Tor^A_2(H_T^*(X,X_0),k)=0$,
% and Lemma~\ref{Tor-vanishing}\,(\ref{Tor-vanishing-2}) together with Lemma~\ref{finitely-generated}
% implies that $\Tor^A_1(H_T^*(X,X_0),k)=0$.
% Hence, the localisation map~$H_T^*(X,X_0)\to S^{-1}H_T^*(X,X_0)$
% is injective by Lemma~\ref{localisation-injective}, and
% Proposition~\ref{localization-iso} gives
% $S^{-1}H_T^*(X,X_0)=S^{-1}H_T^*(X_1,X_0)$.
% The same reasoning as before proves that \eqref{short-sequence} is exact
% for~$i=1$. Continuing this way yields the desired result.
This finishes the proof of Theorem~\ref{main-result}.

\begin{remark}\label{weaker-result}\rm
  Using Remarks \ref{weaker-1}~and~\ref{weaker-2},
  we get a partial result by the same reasoning as above
  if instead of demanding the connectedness of all isotropy groups,
  we make the following weaker assumptions for some~$i_0<n$:
  \begin{enumerate}
  \item $T_x$ is connected for all~$x\in X_{i_0+1}$,
  \item $T_x$ is contained in a subtorus of dimension~$n-(i_0+1)$
    for~$x\not\in X_{i_0+1}$.
  \end{enumerate}
%   that for all~$x\in X\setminus X_i$, $i\le i_0$,
%   the isotropy group~$T_x$ is contained in a subtorus of dimension~$n-(i+1)$.
%   (Note that this assumption in particular means that $T_x$ is connected
%   for all~$x\in X_{i_0}$.)
  Condition~\eqref{Tor-1} % of Theorem~\ref{main-result}
  then implies that the sequences~\eqref{short-sequence} are exact
  for~$i=0$,~\ldots,~$i_0$. This means that we get an exact sequence
  $$
      0
      \longrightarrow H^*_T(X)
      \longrightarrow H^*_T(X_0)
      \longrightarrow H^{*+1}_T(X_1, X_0)
      \longrightarrow \cdots
      \longrightarrow H^{*+i_0}_T(X_{i_0}, X_{i_0-1}).
  $$
\end{remark}

\section{Examples}

\noindent
Here we use integer coefficients again.
As in Remark~\ref{extension-free} we call an $\HBT$-module extended
if it is isomorphic to~$\HBT\otimes_\ZZZ N$ for some $\ZZZ$-module~$N$.

\subsection{Hamiltonian $T$-manifolds}\label{hamiltonian}

\begin{theorem}\label{hamiltonian-extended}
  Let $X$ be a compact Hamiltonian $T$-manifold
  with connected isotropy groups.
  Then $H_T^*(X)$ is an extended $\HBT$-module.

  More precisely, if $F_1$,~\ldots,~$F_s$ are the (finitely many) components
  of~$X^T$, then there are
  even natural numbers~$\lambda_1$,~\ldots,~$\lambda_s$ such that
  $$
    H_T^*(X)=\bigoplus_k \HBT\otimes H^{*-\lambda_k}(F_k).
  $$
\end{theorem}

% Note that over the integers the condition on the fixed point sets
% means that the isotropy group of each non-fixed point is contained
% in a proper subtorus.

The real version of this result (which does not need the
connectivity assumption, \MFcf~Remark~\ref{comment-real})
goes back to Frankel~\cite{Frankel:59},
but see also Atiyah--Bott~\cite{AtiyahBott:84}, Kirwan~\cite{Kirwan:84}
and (for the method of proof used below)
Tolman--Weitsman~\cite{TolmanWeitsman:99}.
Note that our result also works if the cohomology
of the fixed point components has torsion.

\begin{proof}
%   We say that a prime~$p$ is ``bad'' if
%   multiplication by it kills elements in~$H^*(X^T)$
%   or, equivalently,
%   in the extended module~$H_T^*(X^T)=\HBT\otimes H^*(X^T)$.

\def\lopenint#1{\left(#1\right]}
\def\closedint#1{\left[#1\right]}
  Let $\mu\colon X\to\ttt^*$ be a moment map for the $T$-action on~$X$.
  For generic~$\xi\in\ttt$ the function~$f\colon x\mapsto\pair{\mu(x),\xi}$
  is an equivariant Morse--Bott function on~$X$ whose critical set
  is exactly~$X^T$.
% , which has only finitely many components~$F_k$.
  We may assume that the critical values~$c_k=f(F_k)$ are distinct.
  We will prove by induction on~$k$ that $H_T^*(X_k^+)$
  is a finitely generated extended $\HBT$-module,
  where $X_k^+=f^{-1}(\lopenint{-\infty,c_k+\epsilon})$
  with $\epsilon$ so small
  that the interval~$\closedint{c_k-\epsilon,c_k+\epsilon}$
  contains no critical value other than~$c_k$.

  In addition to~$X^+=X_k^+$ introduce
  $X^-=f^{-1}(\lopenint{-\infty,c_k-\epsilon})$,
  the negative disc bundle~$D$ and the negative sphere bundle~$S$ to~$F=F_k$.
  Since the pair~$(X^+,X^-)$ can be retracted to~$(D,S)$, we have
  $$
    H_T^*(X^+,X^-)=H_T^*(D,S)=H_T^{*-\lambda}(F)
  $$
  by the Thom isomorphism, where $\lambda$ is the (even) Morse index of~$F$.
  The composition
  $$
    H_T^{*-\lambda}(F)=H_T^*(X^+,X^-)\to H_T^*(X^+)\to H_T^*(F)
  $$
  is multiplication by the equivariant Euler class~$e_D$
  of the bundle~$D\to F$.
  The latter becomes, after restriction to any point~$x\in F$,
  the product of the non-zero weights of the representation
  of~$T$ on the tangent space to~$D$ at~$x$.
%   In particular, $e_D$ is not divisible by bad primes.
%   Summing up, we get that multiplication by~$e_D$ is injective and does not
%   change divisibility by bad primes.
  In particular, $e_D$ is not divisible by any prime integer,
  so that multiplication by~$e_D$ is injective and does not
  change divisibility by prime integers.

  By injectivity, the long exact sequence for the pair~$(X^+,X^-)$ splits
  into short exact sequences
  $$
    0 \longrightarrow
    H_T^*(X^+,X^-) \longrightarrow
    H_T^*(X^+) \longrightarrow
    H_T^*(X^-) \longrightarrow 0.
  $$
  We claim that this sequence is split. This will imply our claim
  that $H_T^*(X^+)$ is a finitely generated extended $\HBT$-module.

  By induction, $H_T^*(X^-)$ is extended, say of the form~$\HBT\otimes N$,
  where $N$ is the direct sum of some~$H^*(F_i)$, with degree shifts.
  In particular, $N$ is a finitely generated $\ZZZ$-module.
  A section of $\HBT$-modules to~$H_T^*(X^+)\to H_T^*(X^-)$ is the same
  as a section of $\ZZZ$-modules~$N\to H_T^*(X^+)$. The latter exists
  if and only if every element of~$N$ has a preimage
  with the same annihilator.
  % (Of course, this is only a problem for~$R\subsetneq\QQQ$.)
  Take a torsion element~$\gamma\in N$,
  whose annihilator is generated by some~$m\in\ZZZ$.
  % which is a product of bad primes.
  Let $\beta\in H_T^*(X^+)$ be a preimage of~$\gamma$.
  Then $m\beta$ is in the image of some~$\alpha\in H_T^*(X^+,X^-)$.
  Since $m\beta$ is divisible by~$m$, so is $e_D\alpha$,
  hence also~$\alpha$. We finally subtract the image
  of~$\alpha/m$ from~$\beta$ to get the preimage of~$\gamma$
  we were looking for.
\end{proof}

Together with Theorem~\ref{main-result}, this Theorem
%~\ref{hamiltonian-extended}
gives the injectivity result of
Tolman--Weitsman~\cite[Proposition~7.2]{TolmanWeitsman:98}
% see Corollary~\ref{injectivity} below.
% We cannot deduce this from Theorem~\ref{Atiyah-Bredon}
% because $H_T^*(X)$ may not be free (if~$R\subsetneq\QQQ$).
%Together with Corollary~\ref{Chang-Skjelbred-integers},
%Theorem~\ref{hamiltonian-extended} gives
and Schmid's version of the
Chang--Skjelbred lemma \cite[Th\'eor\`eme~3.2.1]{Schmid:01}
in the case of connected isotropy groups. For a discussion
of non-connected isotropy groups see~\cite{FranzPuppe:??};
\MFcf~also~\cite[Section~4]{TolmanWeitsman:99}.

\subsection{Other examples and counterexamples}
\label{examples}

We saw in Remark~\ref{extension-free} that a $T$-space satisfies
condition~\eqref{Tor-j} of Theorem~\ref{main-result}
if its equivariant cohomology is extended.
The converse does not hold in general.

\begin{example}\label{example-2}\rm
  Here we take $T=S^1$.
  Let $X=S^2\cup_\phi C\MFRP^2\times S^1$, where
  $C\MFRP^2$ is the cone over~$\MFRP^2$ and
  $\phi\colon \MFRP^2\times S^1\stackrel{\rm pr}\longrightarrow
  \MFRP^2\to\MFRP^2/\MFRP^1=S^2$.
  The $S^1$-action is trivial on~$S^2$ and free on~$C\MFRP^2\times S^1$
  (acting on the second component), so all isotropy groups are connected.
  It is easy to calculate the (equivariant) cohomology of~$X$
  using the Mayer--Vietoris sequence. One obtains:
  \begin{align*}
    H^i(X) &= \left\{\begin{array}{ll}
        \ZZZ & i=0,2, \\
        \ZZZ/2\ZZZ & i = 4, \\
        0 & \text{otherwise},
      \end{array}\right. \\
    H_T^*(X) &\cong \ZZZ[t]\oplus\ker\bigl(\ZZZ[t]\to\ZZZ/2\ZZZ\bigr)[+2],
  \end{align*}
  where the map $\ZZZ[t]\to\ZZZ/2\ZZZ$ is the composition of the augmentation
  with the projection~$\ZZZ\to\ZZZ/2\ZZZ$ and
  ``$[+2]$'' denotes a degree shift by~$+2$.

  The equivalent conditions \eqref{inclusion-fibre}~to~\eqref{Atiyah-exact}
  are satisfied in this example, but $H_T^*(X)$ is not extended,
  % (in the sense of Remark~\ref{extension-free}),
  in particular, $\iota^*\colon H_T^*(X)\to H^*(X)$ does not have a section.
  % Note that, with rational coefficients,
  % $H_T^*(X;\QQQ)$ is a free $H^*(BT;\QQQ)$-module.\comment{why important?}
\end{example}

In order to get exactness of the Atiyah--Bredon sequence, we assumed that
all isotropy groups are connected.
This requirement cannot be dropped in general.

\begin{example}\rm
  Let~$X=\MFRP^2$, considered as the quotient of the disk~$D^2$
  by identifying opposite points on the boundary~$S^1$.
  The standard rotation of~$T=S^1$ on~$D^2$ descends to~$\MFRP^2$.
  Note that the points on~$\MFRP^1=S^1/\{\pm1\}$ have $\ZZZ_2$ as isotropy group.

  The open sets $\MFRP^2\setminus\MFRP^1$ and $\MFRP^2\setminus0$
  cover~$\MFRP^2$.
  The Mayer--Vietoris sequence for the equivariant cohomology of~$\MFRP^2$
%   =\MFRP^2\setminus\MFRP^1\cup\MFRP^2\setminus0$. 
%   We have $H_T^*(\MFRP^2\setminus\MFRP^1)=\ZZZ[t]$,
%   $H_T^*(\MFRP^2\setminus0)=H_T^*(\MFRP^1)=\ZZZ[t]/2t$
%   and $H_T^*(\MFRP^2\setminus(\MFRP^1\cap0))=H_T^*(S^1)=\ZZZ$.
%   This implies that the Mayer--Vietoris sequence
  splits into short exact sequences, which gives
%   $$
%     0\to H_T^*(\MFRP^2)\to H_T^*()\oplus H_T^*((0,1]\times S^1)
%       \to H_T^*((0,1)\times S^1)\to0.
%   $$
  $$
    H_T(\MFRP^2)=\ZZZ[t]\oplus\ZZZ/2\ZZZ[t][+2]=H^*(\MFRP^2)\otimes \HBT
  $$
  Hence, $H_T^*(\MFRP^2)$ is an extended $\HBT$-module,
  but the map~$H_T^*(X)\to H_T^*(X^T)$ cannot be injective
  because $H_T^*(X^T)=H_T^*(\MFpt)$ is free over~$\ZZZ$.

  By taking products of~$\MFRP^2$ with itself,
  one obtains analogous examples for higher-dimensional~$T$.
\end{example}

The situation is different if one assumes
$H_T^*(X)$ to be free over~$\HBT$.
The main result of~\cite{FranzPuppe:??} is that
in this case one can allow isotropy groups
with at most one (finite) cyclic factor and still get
the exact Atiyah--Bredon sequence.
An example due to Tolman--Weitsman~\cite[Section~4]{TolmanWeitsman:99},
$S^2\times S^2$ with double speed rotation on each factor,
shows that isotropy groups with two cyclic factors are not allowed.
% , even if $H_T^*(X)$ is free over~$\HBT$.

It is not difficult to show that
if $\Tor^{\HBT}_1(H_T^*(X),\ZZZ)=0$ and $H_T^*(X)$ is free over~$\ZZZ$,
then the map~$H_T^*(X)\to H_T^*(X^T)$ is always injective, no matter
what isotropy groups occur. But these assumptions do not guarantee
exactness of the Atiyah--Bredon sequence at~$H_T^*(X^T)$
in the presence of non-connected isotropy groups.

\begin{example}\label{example-4}\rm
  Starting with the $S^1$-space~$X$ from Example~\ref{example-2},
  we consider the space $Y=X\times S^2$ with the componentwise action
  of~$T=S^1\times S^1$, where $S^1$ acts on~$S^2$ with double speed rotation.
  Then $H_T^*(Y)=H_{S^1}(X)\otimes H_{S^1}(S^2)$ is not free over~$\HBT$.
  But since $H_{S^1}(X)$ and $H_{S^1}(S^2)$ are both free over~$\ZZZ$,
  it is easy to check that $H_T^*(Y)\to H_T^*(Y_0)$ is injective.
  Because $H_T^*(Y,Y_0)$ contains an element of order~$4$,
  while $H_T^*(Y_1,Y_0)$ does not,
  $H_T^*(Y,Y_0)\to H_T^*(Y_1,Y_0)$ cannot be injective.
  % (Unlike~$H_T^*(Y_1,Y_0)$, $H_T^*(Y,Y_0)$ contains an element of order~$4$.)
  In other words, the sequence $H_T^*(Y)\to H_T^*(Y_0)\to H_T^{*+1}(Y_1,Y_0)$
  is not exact.
\end{example}

Clearly, for~$T=S^1$ the Atiyah--Bredon sequence is exact
if and only if the restriction map~$H_T^*(X)\to H_T^*(X_0)$ is injective.
Our last example shows
that this equivalence is false for higher-dimensional~$T$,
even if all isotropy groups are connected.

\begin{example}[C.~Allday]\rm
  Let $X=\Sigma T$ be the suspension of~$T$
  with the standard level-wise action. Again the Mayer--Vietoris sequence
  gives the equivariant cohomology:
  \begin{align*}
    \tilde H^*(X) &\cong \tilde H^*(T)[+1], \\
    H_T^*(X) &\cong \ker\Bigl(\ZZZ[t_1,\ldots,t_n]\oplus\ZZZ[t_1,\ldots,t_n]
                              \stackrel{(\epsilon,\epsilon)}\longrightarrow
                              \ZZZ\Bigr),
  \end{align*}
  where $\epsilon\colon\ZZZ[t_1,\ldots,t_n]\to\ZZZ$ denotes the augmentation.
  Note that only for~$n=1$ we have that $H_T^*(X)$ is free over~$\HBT$.
  The map~$H_T^*(X)\to H_T^*(X_0)$ clearly is injective,
  but for~$n>1$ the Atiyah--Bredon sequence is exact only at~$H_T^*(X)$.
\end{example}

\medskip

\footnotesize

\textsc{Fachbereich Mathematik, % und Statistik,
  Universit\"at Konstanz, 78457 Konstanz, Germany}

\emph{E-mail address:}~\texttt{matthias.franz@ujf-grenoble.fr}

\emph{E-mail address:}~\texttt{volker.puppe@uni-konstanz.de}

\end{document}